\documentclass[a4paper,12pt]{amsart}

\usepackage{amsmath}
\usepackage{amssymb,amsbsy,amsmath,amsfonts,amssymb,amscd}
\usepackage{latexsym}
\usepackage{graphics}
\usepackage{color}
\input xy
\xyoption{all}

\newcommand\sC{{\mathcal C}}
\newcommand\sT{{\mathcal T}}
\newcommand\sD{{\mathcal D}}

\newcommand\sF{{\mathcal F}}

\newcommand\sN{{\mathcal N}}

\newcommand\sX{{\mathcal X}}

\newcommand\la{\lambda}
\newcommand\Lam{\Lambda}
\newcommand\al{\alpha}

\newcommand\e{\epsilon}

\newcommand\Ga{\Gamma}

\newcommand\ga{\gamma}
\newcommand\de{\delta}

\DeclareMathOperator{\Mat}{Mat}

\newcommand{\CC}{\ensuremath{\mathbb{C}}}
\newcommand{\RR}{\ensuremath{\mathbb{R}}}
\newcommand{\ZZ}{\ensuremath{\mathbb{Z}}}
\newcommand{\QQ}{\ensuremath{\mathbb{Q}}}

\newcommand{\sS}{\ensuremath{\mathcal{S}}}

\newcommand{\NN}{\ensuremath{\mathbb{N}}}

\newcommand{\PP}{\ensuremath{\mathbb{P}}}

\newcommand{\ra}{\ensuremath{\rightarrow}}

\def\eea{\end{eqnarray*}}
\def\bea{\begin{eqnarray*}}

\DeclareMathOperator{\Id}{Id}

\newcommand{\Proof}{{\it Proof. }}

\newcommand\dual{\mathrel{\raise3pt\hbox{$\underline{\mathrm{\thinspace d
\thinspace}}$}}}
\newcommand\qe{\ifhmode\unskip\nobreak\fi\quad $\Box$}       

\def\BOX{\hfill\lower.5\baselineskip\hbox{$\Box$}}

\newtheorem{theo}[equation]{Theorem}
\newtheorem{remark}[equation]{Remark}
\newenvironment{rem}{\begin{remark}\rm}{\end{remark}}

\newtheorem{defin}[equation]{Definition}

\newtheorem{prop}[equation]{Proposition}

\newtheorem{example}[equation]{Example}

\DeclareMathOperator{\Irr}{Irr}

\DeclareMathOperator{\GL}{GL}

\begin{document}

\title[Moduli of GH Manifolds]{Teichm\"uller  spaces of Generalized Hyperelliptic manifolds}
\author{ Fabrizio  Catanese - Pietro Corvaja}
\address {Lehrstuhl Mathematik VIII\\
Mathematisches Institut der Universit\"at Bayreuth\\
NW II,  Universit\"atsstr. 30\\
95447 Bayreuth}
\email{fabrizio.catanese@uni-bayreuth.de}
\address {Dipartimento di Scienze Matematiche, Informatiche e Fisiche\\
Universit\'a di Udine\\
via delle Scienze, 206\\
33100 Udine}
\email{ pietro.corvaja@uniud.it}

\thanks{
The present work took place in the framework  of the 
 ERC Advanced grant n. 340258, `TADMICAMT' }

\date{\today}
\maketitle

\begin{abstract}
In this paper we answer two questions posed by  \cite{topmethods}, thus achieving in particular  a description of 
the connected components of Teichm\"uller space corresponding  to  Generalized Hyperelliptic Manifolds $X$.  These are the quotients $ X = T/G$
of a complex torus $T$ by the free action of a finite group $G$,
and they are also the K\"ahler classifying spaces for a certain class of Euclidean cristallographic groups $\Ga$,
the  ones which are torsion free and even.
\end{abstract}

\tableofcontents

\section*{Introduction}

The classical hyperelliptic surfaces  are the quotients of a complex torus of dimension 2  by a finite 
group $G$ acting freely, and in such a way that the quotient is not again a complex  torus.

These surfaces were classified by Bagnera and de Franchis (\cite{BdF}, see also \cite{Enr-Sev} and \cite{bpv}) and they are obtained as quotients $(E_1 \times E_2)/G$
where  $E_1, E_2$ are  two elliptic curves, and $G$ is an abelian group acting on $E_1$ by translations, and on $E_2$  effectively and  in such a way that $ E_2/G \cong \PP^1$.

In higher dimension we define the 
Generalized Hyperelliptic Manifolds (GHM) as quotients $T/G$ of a (compact) complex torus  $T$ by a finite group $G$ acting freely,
and with the property  that $G$ is not a subgroup of the group of  translations. Without loss of generality one can then assume
that $G$ contains no translations (since the subgroup $G_T$ of translations in $G$ would be a normal subgroup, and
if we denote $G' = G/G_T$, then $T/G = T' / G'$, where $T'$ is the torus  $ T' : = T/G_T$).

The name  Bagnera-de Franchis (BdF) Manifolds  was instead reserved in \cite{topmethods} and \cite{bcf} for those quotients $X = T/G$ were $G$ contains no translations,
and $G$ is a cyclic group of order $m$ (observe that, when $T$ has  dimension $n=2$, the two notions coincide, thanks
to the classification result of Bagnera-De Franchis in \cite{BdF}). BdF manifolds of small dimension were studied in \cite{topmethods} and \cite{demleitner}.

Before stating our  main theorem,
recall first of all  that the Teichm\"uller space $\sT (X)$ of a compact complex manifold $X$  is the quotient
$$\sT (X) : = \sC 
\sS (X) / \sD iff^0 (X)$$
of the space of complex structures on the  oriented differentiable manifold underlying $X$, which are compatible with the natural orientation
of $X$, by the diffeomorphisms of $X$ which are isotopic to the identity.

 Recall also that  a $K(\Ga, 1)$ manifold is a manifold $M$ such that its universal covering is contractible, and such that $\pi_1(M) \cong \Ga$.
 
 We have then:

\begin{theo}\label{mt}
Given a Generalized Hyperelliptic Manifold $X$, $X$ is K\"ahler and  its fundamental group $\pi_1(X)$  is a torsion free  even Euclidean cristallographic group $\Ga$ (see definitions \ref{cristall} and  \ref{even}).

Conversely, given such a torsion free  even Euclidean cristallographic group $\Ga$, there are GHM with $\pi_1(X) \cong \Ga$; moreover 
 any compact K\"ahler manifold $X$ which is a $K(\Ga, 1)$ is a Generalized Hyperelliptic manifold.

The subspace of the Teichm\"uller space $\sT (X) $ corresponding to K\"ahler manifolds consists of a  finite number  of connected components, indexed  by the Hodge type
of the Hodge decomposition. Each such component is a product of open sets of complex Grassmannians.
\end{theo}

\section{Euclidean cristallographic groups}

\begin{defin}\label{cristall}
(i) We shall say that a group $\Ga$ is an abstract Euclidean cristallographic group if there exists an exact
sequence of groups
$$ (*) \ \ 0 \ra \Lam \ra \Ga \ra G \ra 1 $$
such that
\begin{enumerate}
\item
$G$ is a  finite group
\item
$\Lam$ is free abelian (we shall denote  its rank by $r$)
\item
Inner conjugation $ Ad : \Ga \ra  Aut (\Lam) $ has Kernel exactly $\Lam$,
hence $ Ad$  induces an embedding, called {\bf Linear part},
$$ L : G \ra GL (\Lam) : = Aut (\Lam) $$
(thus $L(g) (\la) =  Ad (\ga) (\la) = \ga \la \ga^{-1}, \ \forall \ga {\rm \ a \ lift \ of }\  g$) 

\end{enumerate}

(ii) An {\bf affine realization defined over a field $ K \supset \ZZ$ }  of an abstract Euclidean cristallographic group $\Ga$  is a homomorphism
(necessarily injective) 
$$\rho : \Ga \ra Aff (\Lam \otimes_{\ZZ} K)$$ such that 

[1] $\Lam$ acts by translations on $ V_K := \Lam \otimes_{\ZZ} K$,  $ \rho(\la) (v) =   v + \la$,

[2]  for any $\ga$ a lift of $g \in G$ we have:
$$  V_K \ni v \mapsto \rho(\ga) (v) =  Ad (\ga) v + u_{\ga} = L(g) v + u_{\ga}, \ {\rm for \ some } \ \ u_{\ga} \in  V_K.$$ 

(iii) More generally we can say that an affine realization of $\Ga$  is obtained via a lattice $\Lam ' \subset \Lam \otimes_{\ZZ} \QQ$
if there exists a homomorphism $\rho ' : \Ga \ra Aff (\Lam ')$ such that $\rho = \rho ' \otimes_{\ZZ}  K$ (then necessarily $\Lam \subset \Lam'$).

\end{defin}

\begin{rem}\label{notation}
In the previous formulae in (ii) [2] and in the following we used a shorthand notation, we extend the action of  $L(g)$ on $\Lam$
to $V_K$ naturally as $ L(g)\otimes_{\ZZ} Id_K$. We shall also often write 
$ g (v) : = L(g) (v)$, and $\ga (v) = Ad (\ga) (v)$.

 We note that for a cristallographic group $\Gamma$,  realizing it via a lattice $\Lam '$  as in (iii) amounts to having all the $u_\gamma$ inside the lattice $\Lam '$, in the formula appearing in (ii) [2].
\end{rem}

\begin{rem}\label{unique}
Given  a Euclidean cristallographic group $\Ga$ as above, the exact sequence $(*)$  is unique up to isomorphism, since $\Lam$ is the 
unique maximal normal abelian subgroup 
of $\Ga$ of finite index. 

In fact, if $\Lam'$ has the same property,  then their intersection $\Lam^0 : = \Lam \cap \Lam'$ is a normal subgroup of
 finite index, in particular $\Lam^0  \otimes_{\ZZ} \QQ = \Lam  \otimes_{\ZZ} \QQ = V_{\QQ} $
  and any automorphism of $\Lam$ which is the identity on $\Lam^0$ is the identity.
 
Since $\Lam' \subset \ker ( Ad  : \Ga \ra \GL  (\Lam^0))$, $\Lam' \subset \ker ( Ad  : \Ga \ra \GL  (\Lam))  =  \Lam$:
by maximality $\Lam' = \Lam$.

\end{rem}

\begin{rem}\label{cocycle}
 $L$ makes $\Lam$ and $V_K$  left $G$-modules, and to give an affine realization is equivalent to giving 
 a 1-cocycle in 
$Z^1 ( \Ga , V_K)$ such that $u_{\la} = \la \ \forall \la \in \Lam$, since 
$$ \rho (g_1 g_2) =   \rho (g_1 )  \rho ( g_2)  \iff  L(g_1 ) ( L(g_2) v + u_{\ga_2})  + u_{\ga_1} =  L(g_1 g_2) v + u_{\ga_1 \ga_2} \iff$$
$$ \iff     u_{\ga_1 \ga_2} =  u_{\ga_1}  +   L(g_1 ) (  u_{\ga_2}) =   u_{\ga_1}  +   g_1  (  u_{\ga_2}) .$$ 

(c) Two such cocycles $(u_{\ga}), (u'_{\ga}),$ are cohomologous if and only if there exists a vector $w \in V_K$ such that:
$$  u'_{\ga} - u_{\ga} = \ga w - w, \ \forall \ga \in \Ga.$$ 

(d)  Hence two such  cocycles are cohomologous  if and only if the respective affine realizations are conjugate by a translation in $V_K$, since
$$ \rho(\ga) (v + w) - w = \ga (v + w)  + u_{\ga} - w =  \ga v  +  u_{\ga} + (\ga w  - w) =  \ga v  +  u'_{\ga} =  \rho' (\ga) (v).$$ 

\end{rem}

\begin{theo}\label{affinereal}
Given an abstract Euclidean cristallographic group there is a unique class of affine realization, for each    field $K \supset \ZZ$.

There is moreover an effectively computable  minimal number $d \in \NN$ such that the  realization is obtained via $\frac{1}{d} \Lam$.
\end{theo}

\Proof
$Ad : \Ga \ra GL(\Lam)$ makes $\Lam$ a $\Ga$-module, a trivial $\Lam$ module, hence also a $G$-module.

We have seen in remark \ref{cocycle}, (b), that an affine realization is given by a cocycle $u_ {\ga} $ in $Z^1(\Ga, V_K)$
such that $u_{\la} = \la, \ \forall \la \in \Lam$; and moreover the class of the realization depends only on the cohomology class
in   $H^1(\Ga, V_K) $.

 Consider now the exact sequence of  cohomology groups
$$ H^1(G, V_K) \ra H ^1(\Ga, V_K) \ra H^1(\Lam , V_K) = Hom (\Lam , V_K) \ra H^2(G, V_K) .$$

Since $G$ is a finite group and $K$ is field of characteristic zero,   $ H^1(G, V_K) = H^2(G, V_K)  = 0$ (\cite{jacobson} pages 355-363): hence 
we get an isomorphism $H ^1(\Ga, V_K) \ra H^1(\Lam , V_K) = Hom (\Lam , V_K) $.

We look for a cohomology class  $[u_ {\ga} ]$ such that its image in  $H^1(\Lam , V_K) = Hom (\Lam , V_K)  $
is the tautological map $ \la \mapsto \la \in V_K$, composition of
the  identity of $\Lam$ with the inclusion $\Lam \subset V_K$. By the above isomorphism 
such  cohomology
class exists,  is unique and not equal to zero.

In particular, this applies for $ K = \QQ$, and  since $G$ is finite there is an integer $D$ such that 
that  $u_ {\ga}  \in \frac{1}{D} \Lam $. 

Hence we get a cocycle  in $ H ^1(\Ga, \frac{1}{D}   \Lam) $  whose image in $ H ^1(\Ga, (\frac{1}{D} \Lam) / \Lam ) $
comes form a unique class in $ H^1(G, (\frac{1}{D} \Lam) / \Lam ) $.
This last class has order dividing $D$, and we let $\de$ be its order. Hence,
setting $ D = d \de$ we obtain that  $[u_ {\ga}]$ comes from  $  H ^1(\Ga, (\frac{1}{d} \Lam))$.

\qed

\begin{rem}\label{erratum}
  The above result provides a correct proof for a claim made in  \cite{bc-inoue} and \cite{bcf},    answers the   question posed in \cite{topmethods},
remark 23 , page 312, and  generalizes the result of  lemma 18, page 310 ibidem.

We realized later on  that the first statement, about the unicity of the affine realization, was proven by Bieberbach in 1912 (\cite{bieb2}).
 \end{rem}
 
  We have moreover
\begin{prop}\label{splitting}
 Suppose that $\Gamma$ is an Euclidean cristallographic group, which by theorem  \ref{affinereal} is a subgroup of $Aff(V_\QQ)$.   Consider a lattice $\Lam ' \subset \Lam \otimes_{\ZZ} \QQ$,
 such that $\Lam \subset \Lam'$,  and let $\Ga'$ be the corresponding group of affine transformations of $V_\QQ$, generated by $\Ga$ and $\Lam '$. Then the group $\Ga'$ sits into a canonical exact sequence 
 $$0 \ra \Lam'  \ra \Ga' \ra G \ra 1 ,$$
and such an exact sequence splits if and only if the affine realization of $\Ga$    is defined over $\Lam'$.
 \end{prop}

\Proof
 The group $\Ga'$ is obtained from $\Ga$ by adding some translations, which lie in the kernel of $Ad: Aff(V_\QQ)\to GL(V_\QQ)$, so the image $Ad(\Ga')$ coincides with $G=Ad(\Ga)$.  Hence we get the exact sequence above.

 Assume now that the affine realization is defined over $\Lam'$. Then, for each $g\in G$, let $\ga\in \Ga \subset \Ga'$ be any lift of $g$. From 
$\gamma (v) = L(g) v + u_\gamma$ for some $u_\gamma \in\Lam '$, we obtain a second lift $\gamma'\in \Ga'$ of $g$ by setting $\gamma'= u_\gamma^{-1}  \circ \ga $, i.e. $\gamma' (v)= L(g) v$. Clearly the homomorphism $g \mapsto \ga'$ yields a splitting.

Conversely, if the exact sequence splits, then we have a semidirect product of $\Lam'$ with a group $G'$ isomorphic to $G$. Then there is a fixed point $w$ for the action
of $G'$ on $V_\QQ$, obtained as $w : = \Sigma_{g \in G'} gv$, where $v$ is any vector in $\Lam$.
Choosing $w$ as the origin, we obtain an   affine realization of $\Ga'$ defined over $\Lam'$, a fortiori the same holds for $\Ga$. 

\qed

\begin{prop}\label{prop}

(I) Let $\e \in H^2 (G, \Lam)$ be the extension class of $\Ga$. 

Then
 there is an affine realization of $\Ga$ defined over $\Lam'$ 

$\Leftrightarrow$ $\e \otimes_{\ZZ}  \Lam ' = 0$ 

$\Leftrightarrow$  there is a fixed point $ [w] \in \Lam' / \Lam $ for the action of $G $ on $V_\QQ / \Lam$.

(II) $G$ acts freely on the real torus $T : = ( \Lam \otimes_{\ZZ} \RR ) / \Lam = V_{\RR} / \Lam $ if and only if $\Ga$ is {\bf torsion free}: this means  that  the subset
$$ Tors(\Ga) : = \{ \ga | \exists m \in \NN_+ , \ga^m = 1_{\Ga} \} $$ consists only of the identity.

\end{prop}

\Proof

(I) An extension splits if and only if its cohomology class $\e'  \in H^2 (G, \Lam ') = 0$ (\cite{jacobson} theorem 6.15, page 365),
 hence the first assertion follows from proposition \ref{splitting}.

The second assertion follows as in the proof of   proposition \ref{splitting}.

(II) if $g$ acts with a fixed point on $T$, there is a lift $\ga$ of $g$ such that $\ga$ has a fixed point $w $  in $V$. Then, choosing 
coordinates such that $ w =0$, the action of $\ga$ is linear, hence the order of $\ga$ equals the order of $g$.

Conversely, if $\ga$ has finite order $m$, the vector $ w : = \sum_1^m \ga^i (0)$ is fixed by $\ga$, hence $[w]$ is fixed by $g$.

\qed
 
 \begin{rem}\label{BdF}
 Let $X = T/G$  be a Generalized Hyperelliptic Manifold.
The action of   $g \in G$  is induced by an affine transformation $ x \mapsto \al x + b $ on the universal cover, hence it 
does not have a fixed point on $ T =  V/ \Lambda$ if and only if there is no solution of the equation
$$   g (x) \equiv  x  \  \pmod  \Lambda $$ 
to be solved in $x\in V$, i.e.   to the equation
$$  \la \in \Lambda, \  (\al - \Id) x = \la - b$$
in $(x,\lambda)\in V\times \Lambda$.
This remark shows that if the action of $G$ on $X$ is free it is necessary that  $1$   be an eigenvalue of $\al=L(g)$ for all non trivial transformations $ g \in G$.
 \end{rem}

\section{Actions of a finite group on a complex torus $T$}

Assume that we have the  action of a finite group $G'$ on  a complex torus 
 $ T= V / \Lam'$, where $V$ is a complex vector space, and $ \Lam' \otimes_{\ZZ} \RR \cong V$.
 
   Since every holomorphic map between complex tori lifts to a complex  affine map
 of the respective universal covers, we can attach  to the group $G'$  the group $\Ga$  of (complex) affine transformations
 of $V$ which are  lifts of  transformations of the group $G'$.

 \begin{prop}\label{notranslations}
 $\Ga$ is an Euclidean cristallographic group, via the  exact sequence 
 $$0 \ra \Lam \ra \Ga \ra G \ra 1  $$
where  $\Lam  \supset  \Lam'$  is the  lattice in $V$ such that  $\Lam: = Ker (Ad), \ Ad : \Ga  \ra GL  (\Lam')$.

Defining  $G^0$ to be the subgroup of $G'$ consisting of all the translations in $G'$,  then  $  G^0 = \Lam / \Lam' $,
and moreover $ G \subset Aut (V/ \Lam)$
 contains no translations.
  \end{prop}
 
 \Proof
$\Lam$ is a subgroup of the group of translations in $V$, hence it is obviously Abelian,  and maps isomorphically onto a lattice of $V$ 
which contains $\Lambda'$.  We shall identify this lattice with $\Lambda$.

 $\Lam$  is normal, $  G^0 = \Lam / \Lam' $  and $G $ embeds in $ GL (\Lam') \subset GL (\Lam) $. 
 
\qed

Hence the datum of  the action of a finite group  $G'$ on a complex torus $T$ is equivalent to giving:

\begin{enumerate}
\item
a cristallographic group $\Ga$
\item
 a $G$- invariant  sublattice $\Lam'$ of the maximal normal abelian subgroup $\Lam$ of finite index
 (equivalently, to give a normal such sublattice $\Lam'$), so that we may set   $G' : = \Ga / \Lam'$
\item
 a complex structure $J$ on the real vector space $V_{\RR}$ which makes the action of $G$ complex  linear. 
\end{enumerate}

While the data (1) and (2) are discrete data, the choice of the complex structure $J$ on $V$ may give rise to 
positive dimensional  moduli spaces. We introduce  however (see \cite{topmethods}) a further discrete  invariant, called Hodge type.

\begin{defin}\label{HodgeT}
(i) Given a  faithful representation $ G \ra Aut (\Lam)$, where $\Lam$ is a free abelian group  of even rank $r = 2n$,
 a {\bf $G$- Hodge decomposition} is a $G$-invariant decomposition 
 $$\Lam \otimes \CC = H^{1,0} \oplus H^{0,1}, \ H^{0,1} = \overline{H^{1,0}}.$$

(ii)
Write   $\Lam \otimes \CC$
as the sum of
isotypical components $$\Lam \otimes \CC = \oplus_{\chi \in \Irr (G)}
U_{ \chi}.$$
Write also  $U_{ \chi} = W_{ \chi} \otimes M_{ \chi}$, where $W_{ \chi} $ is the  irreducible representation corresponding to the character $\chi$,
and $ M_{ \chi} $ is a trivial representation whose  dimension is denoted $n_{ \chi} $.

Write accordingly $V : = H^{1,0} =  \oplus_{\chi \in \Irr (G)}
V_{ \chi},$ where $V_{ \chi} = W_{ \chi} \otimes M^{1,0}_{ \chi}$.

 Then the {\bf Hodge type} of the decomposition
is the datum 
of the dimensions   $$\nu ( \chi): = dim_{\CC} M^{1,0}_{ \chi} $$
corresponding to  the Hodge summands for non real representations (observe in fact  that one must have: $\nu ( \chi) + \nu ( \bar{\chi}) = \dim ( M_{ \chi})$).
\end{defin}

\begin{defin}\label{even}
A cristallographic group $\Ga$ is said to be {\bf even} if:

i)  $\Lam$ is a free abelian group  of even rank $ r = 2n$

ii) considering the associated   faithful representation $ G \ra Aut (\Lam)$,  for each real representation $\chi$, 
 {\bf $ M_{ \chi}$ has even dimension } (over $\CC$).

\end{defin}

\begin{rem}\label{Hodgetype}
(i) Given a group action on a complex torus, we obtain an Euclidean  cristallographic group which is even, since $\Lam \otimes_{\ZZ} \RR$
admits a $G$- invariant complex structure.

 (ii) Given an even cristallographic group $\Ga$,  
the $G$- Hodge decompositions of a fixed Hodge type (satisfying the necessary condition
$\nu ( \chi) + \nu ( \bar{\chi}) = \dim ( M_{ \chi})$) are parametrized by a (non empty)  open set in a product of Grassmannians, as follows.
 
 For a non real irreducible representation $\chi$ one may simply choose $M^{1,0}_{ \chi}$ to be a complex subspace of dimension  $\nu ( \chi)$ of  $ M_{ \chi}$, and for $M_{ \chi} = \overline{M_{ \chi}}$,
one simply chooses a complex subspace $M^{1,0}_{ \chi}$ of middle dimension. 

They must satisfy the open condition  that (since $  M^{0,1}_{ \chi} : = \overline {M^{1,0}_{ \chi} }$) 
$$ M_{ \chi} = M^{1,0}_{ \chi}  \oplus  M^{0,1}_{ \chi}  \iff M_{ \chi} = M^{1,0}_{ \chi} \oplus  \overline{M^{1,0}_{\bar{ \chi}} } .$$ 
\end{rem}

 \section{Proof of theorem \ref{mt}}
 
1)  Let $X = T/G$ be a GHM.  Then $X$ is K\"ahler, since,  averaging (by the action of $G$) a K\"ahler metric on $T$ with constant coefficients,
 we obtain a $G$-invariant one, which descends to the quotient manifold $X$.
 
 2) Since the universal covering of $X$ is the vector space $V$, which is contractible,
 and $ X = V / \pi_1 (X) $, where $\pi_1 (X) $ acts freely, we obtain, setting $\pi_1 (X) = :  \Ga$  that $X$ is a $ K ( \Ga,1)$ manifold.
 
 3) That  $\Ga$ is an Euclidean cristallographic group follows from proposition \ref{notranslations}, 
  that $\Ga$ is torsion free by (II) of proposition \ref{prop}, and $\Ga$ is even by (i) of  remark \ref{Hodgetype}.
  
   4) If $X'$ is also  a compact K\"ahler manifold which is a  $ K ( \Ga,1)$ , then $X'$ admits a Galois unramified covering $T'$ with group $G$ such that
 $T'$ is a compact K\"ahler manifold with the same integral cohomology algebra of a complex torus, hence $T'$ is a complex torus,
 as shown in \cite{catAV}, see also \cite{topmethods}. Hence also $X' = T' / G$ is a GHM.

  5) Moreover, given an Euclidean cristallographic group such that $\Ga$ is torsion free  and even, $V_{\RR}$ admits a complex structure
  by (i) and (ii) of remark \ref{Hodgetype}, and the action of $\Ga$ on $V$ is free since $\Ga$ is torsion free. Hence for any such complex structure
  we obtain a quotient $ X = V / \Ga = (V/\Lam)/ G = T/G$ which 
  is a GHM. 
  
  6) Observe that   the family of GHM is stable by deformation in the large. Indeed, every deformation of a GHM $X = T/G$  yields a deformation of the covering torus $T$
   (observe that the covering $ T \ra X$ is associated to the unique surjection $\pi_1 (X) \cong \Ga \ra G$)
 and that in  \cite{as}, \cite{nankai} and \cite{cat04},  was proven that a deformation in the large of a complex torus is a complex torus,
 so that the natural family of $n$-dimensional complex tori is a connected component $\sT_n$ of the Teichm\"uller space $\sT (T)$
 (but not the only one). 
 
 7)  The  connected component $\sT_n$ of the Teichm\"uller space of  $n$-dimensional complex tori  (see  \cite{nankai}, 
\cite{cat04} and \cite{handbook} )
is  the open set $\sT_n$ of the complex
Grassmann Manifold $Gr(n,2n)$, image of the open set of matrices

$\sF : = \{ \Omega \in \Mat(2n,n; \CC) \ | \ i^n det  (\Omega \overline
{\Omega}) > 0 \}.$

Over $\sF$ lies the following tautological family of complex tori:    consider
a fixed lattice $ \Lam  : =  \ZZ^{2n}$, and associate
to each matrix $ \Omega $ as above   the subspace $V$ of $\CC^{2n} \cong \Lam \otimes \CC$ given as
$$ V : =  \Omega \CC^{n},$$ so that
$ V \in Gr(n,2n)$ and $\Lam \otimes \CC \cong V \oplus \bar{V}.$

To  $ V $ we associate then the  torus 
$$T_V : = V / p_V (\Lam) = (\Lam \otimes \CC)/ (\Lam \oplus  \bar{V} ),$$
$p_V : V \oplus
\bar{V} \ra V$ being the projection onto the first summand.

The crystallographic group $\Ga$ determines an action of $G \subset SL (2n, \ZZ)$ on $\sF$ and on $\sT_n$,
obtained by multiplying the matrix  $\Omega$ with matrices $g \in G$ on the right.

Define then $\sT_n^G$ as the locus of fixed points for the action of $G$. If $V \in \sT_n^G$, then 
$G$ acts as a group of biholomorphisms of $T_V$, and we associate then to such a $V$ the GHM
$$ X_V : =  T_V / G.$$

Since the induced family $\sX \ra \sT_n^G$ is differentially trivial, we obtain a map
$\psi : \sT_n^G \ra \sT (X)$.

8)  We see that  $\sT_n^G$ consists of  a finite number of  components, indexed by the Hodge type of the Hodge decomposition.
Observe in fact that the Hodge type
 is invariant by deformation, so it distinguishes  a finite number of connected components of $\sT_n^G$.
That these connected components are just   a product of Grassmannians follows from (ii) of remark \ref{Hodgetype}.
 
9)  Assume in greater generality
 that we have an unramified Galois cover $ p : Y \ra X$ such that the associated subgroup $\pi_1(Y) = : \Lam$
is a characteristic subgroup of $\pi_1(X) = : \Ga$, and denote by $G$ the quotient group.
Then, via pull back, the space $\sC \sS (X)$ of complex structures $J'$ on $X$ is contained in the space 
$\sC \sS  (Y)$ of complex structures $J$ on $Y$, and is actually equal to the fixed locus of $G$,
$$ \sC \sS  (X) = \sC \sS  (Y)^G = \{ J | g_* (J) = J , \ \forall g \in G\} = \{ J | G \subset Bihol (Y_J)\}.$$ 

10) Since $\Lam$ is a characteristic subgroup, all diffeomorphisms of $X$ lift to $Y$, and we have an exact sequence

$$  1 \ra G \ra \sN_Y(G) \ra \sD iff (X) \ra 1, \   \sN_Y(G): = \{ \phi \in \sD iff (Y) |  \phi G = G \phi \},$$
since the diffeomorphisms in the normalizer $ \sN_Y(G)$ of $G$ are the deffeomorphisms  which descend to $X$.

11) If $X,Y$ are classifying spaces,  then $\sD iff (X)^0$ is the subgroup acting trivially on $\pi_1(X) =  \Ga$, 
and similarly for $Y$.

12) In our case $G$ acts non-trivially on the first homology, hence we get an inclusion 
$$ \sD iff (X)^0 \subset \sD iff (T)^0,$$
as the normalizer subgroup $\sN_Y(G)^0 $ of $G$. 

13) We consider now Teichm\"uller space $\sT (X) = \sC \sS (X) / \sD iff (X)^0$.
Because of 9) and 12), 
$$\sT (X) = \sC \sS (X) / \sD iff (X)^0 =     \sC \sS  (Y)^G /  sN_T(G)^0 . $$
We get therefore a continuous map $ j : \sT (X) \ra  \sT (T)^G$, where $\sT (Y)^G$ is the image of
$ \sC \sS (Y)^G$ inside $\sT (Y) $.

14) We want to show that in our case  $j$ is a homeomorphism, at least when restricted to the inverse image of $\sT_n$,
which we denote by $\sT (X)_{GH}$.

It suffices to observe  that $j$ and $\psi$ are inverse to each other, and to show that they are  local homeomorphisms.
We use remark 14 and proposition 15 of \cite{handbook}. In fact, locally for the torus $T$ such that $X = T/G$,
Teichm\"uller space $\sT (T) $ is locally the Kuranishi space $ Def (T)$, and in turn $Def(X)$
is the closed subspace $Def(X) = Def(Y)^G$ of fixed points for the action of $G$.

Locally there is a surjection $Def (X) \ra \sT (X) $. Composing it with $j$ we get the composition of the inclusion
$Def(X) \subset  Def(T)$ with the local homeomorphism $Def(T) \cong \sT (T) $: hence the composition is injective.

By 3) remark 14 of \cite{handbook} $Def (X) \ra \sT (X) $ is a local homeomorphism and we can also conclude that $j$ is a homeomorphism with its
image $\sT (T)^G$.

 \bigskip

\subsection{Concluding remark}
We raise the  following question: given a generalized Hyperelliptic Manifold $X$, classify the projective manifolds which are a deformation of $X$
(see \cite{demleitner} for the special case of Bagnera de Franchis manifolds and \cite{U-Y}and \cite{Lange} 
for results in dimension $3$).


\end{document}